\documentclass[12pt]{article}

\usepackage{amsmath, amssymb}

\usepackage{hyperref}

\hypersetup{
colorlinks=true, 
breaklinks=true, 
urlcolor= blue, 
linkcolor= blue, 
bookmarksopen=true, 
pdftitle={}, 
pdfauthor={Christian Bourdarias}, 
pdfsubject={} 
}

\newtheorem{theorem}{Theorem}
\newtheorem{proposition}{Proposition}
\newtheorem{definition}[theorem]{Definition}
\newtheorem{remark}{Remark}

\newcommand{\pr}{{\bf \textit{Proof: }}}
\newcommand{\cqfd}{{\nobreak\hfil\penalty50\hskip2em\hbox{}\nobreak\hfil
$\square$\qquad\parfillskip=0pt\finalhyphendemerits=0\par\medskip}}

\title
{    
      Oscillating  
      waves  and
    \\     optimal  smoothing  effect 
\\ for  one-dimensional  
\\ nonlinear scalar conservation laws }

 \author
{Pierre Castelli and  St\'{e}phane Junca}

\begin{document}

  \date{\today}

 \maketitle

\centerline{\scshape Pierre Castelli }
\medskip
{\footnotesize
 \centerline{Lyc\'ee Jacques Audiberti}
 \centerline{63 boulevard Wilson}
   \centerline{ 06600, Antibes, France}
} 

\medskip

\centerline{\scshape St\'ephane Junca}
\medskip
{\footnotesize
 \centerline{Labo. J. A. Dieudonn\'e, UMR CNRS  7351 }
   \centerline{Universit\'e de Nice Sophia-Antipolis} 
   \centerline{Parc Valrose,06108 Nice cedex 02, France}
\centerline{and}
\centerline{Team Coffee}
\centerline{ INRIA Sophia Antipolis Mediterran\'ee}
\centerline{2004 route des Lucioles - BP 93}
\centerline{06902 Sophia Antipolis Cedex, France}
}

\bigskip

 \centerline{(Communicated by the associate editor name)}

\begin{abstract} 
\bigskip
Lions, Perthame, Tadmor  conjectured  in 1994
an optimal  smoothing effect for entropy solutions of nonlinear scalar conservations laws 
(\cite{LPT}). 
In this short paper  we will  restrict our attention to  the simpler one-dimensional case.  
First,  supercritical geometric optics    lead to  sequences of $C^\infty$ solutions uniformly bounded in the Sobolev space conjectured. 
Second we give  continuous  solutions which belong exactly to 
  the  suitable Sobolev space.
 In order to do so  we give two new definitions of nonlinear flux and we introduce  fractional $BV$ spaces.
\end{abstract}

\medskip


\tableofcontents  

\section{Introduction and nonlinear flux  definitions 
             }\label{sI}
\smallskip
$\mbox{ }$\\

   We focus on oscillating  smooth solutions  for one-dimensional scalar conservations laws:  
\begin{equation} \label{eqcl}
\frac{\partial   u}{\partial t } + \frac{\partial  f( u)}{\partial x } = 0, \qquad u(0,x)=u_0(x), \qquad  t>0, \, x \in \mathbb{R}.
\end{equation}
 The aim  of this paper is to build  solutions  related to  the maximal regularity
 or  the uniform Sobolev bounds  conjectured in \cite{LPT}  for entropy solutions.
In the one-dimensional case,  piecewise smooth solutions 
with the maximal regularity are obtained in \cite{DW}  for power-law fluxes. 
We  seek    supercritical  geometric optics expansions
and  some special oscillating solutions. 
Our constructions are valid  for all   $C^\infty$ flux
 and show  that  one cannot expect a better smoothing effect. 
\\

The more complex multidimensional case is dealt with  in  
\cite{Ju4, CJ2}.   For  recent other approaches we refer the reader to \cite{Da85,Ch00, DOW, COW,JaX09,GP}. 
Recall  that the first famous  $BV$ smoothing effect for uniformly convex flux was given by the  Oleinik one-sided Lipschitz condition 
in  the 1950s  (see for instance the  books \cite{Da,Laxbook}).  
For solutions with bounded entropy production, the smoothing effect is weaker than for entropy solutions (\cite{DW,GP}).
\\

Let us give various definitions of nonlinear flux from \cite{LPT,BJ,Ju4,BGJ6}.
Throughout the paper,   $K$ denotes  a compact  real interval.
\begin{definition}{\bf [Lions-Perthame-Tadmor nonlinear flux,\cite{LPT}]} \label{defLPT}
$\mbox{ }$ \\
 $f \in C^1(K,\mathbb{R})$ is said to be a nonlinear flux on $K$ with degeneracy $\alpha$
 if there exists a constant $C>0$   such that  for all $\delta >0$, 
\begin{equation}\label{nfLPT}
   \sup_{\tau^2+\xi^2=1} \left(   \mbox{measure} \{ v \in K, \;  |\tau + \xi \, f'(v)| < \delta   \} \right)  \leq C  \delta^\alpha.      
\end{equation}
\end{definition}
In \cite{LPT}, the authors proved a  smoothing effect  for entropy solutions
in  some Sobolev space. 
  They obtained  uniform Sobolev bounds  
with respect to  $L^\infty$ bounds of   initial data. 
Moreover, they conjectured a better smoothing effect :
\begin{equation} \label{conj}
 u_0 \in L^\infty(\mathbb{R})   \Rightarrow  u(t,.) \in W^{s,1}_{loc}(\mathbb{R}_x), \; \mbox{ for  all} \;    s < \alpha   \mbox{ and  for  all}\; t> 0 
\end{equation}
where the parameter $\alpha$ is defined in Definition \ref{defLPT}.
They proved a weaker smoothing  effect  which was improved in \cite{TT}.
The conjecture (\ref{conj})  is still an open problem.\\%

In \cite{Ju4} was given   another  definition 
related to the derivatives of the flux. 
It generalizes a notion of nonlinear flux arising in  geometric optics  (\cite{CJR}).
 The   next one-dimensional definition of smooth nonlinear flux is simpler than 
 in the multidimensional case (\cite{BJ,Ju4}). 
\begin{definition}{\bf[Smooth nonlinear flux, \cite{Ju4}]} \label{defJ4}
$\mbox{ }$ \\
 $f \in C^\infty(K,\mathbb{R})$ is said to be a nonlinear flux  on $K$ with degeneracy $d$
 if 
\begin{equation}\label{eqJ4}
  d= \max_{u \in K} \left(  \min \left \{   {k\geq 1}, \;  \frac{d^{1+k} f}{d u^{1+k}} (u) \neq 0  \right \}  \right)   < + \infty.
\end{equation}
\end{definition}
For the Burgers equation or  for uniformly convex flux, the degeneracy is  $d=1$. 
That is the minimal  possible value. 
 For the cubic flux $f(u)=u^3$ on $K=[-1,1]$, the degeneracy is  $d=2$. 
The cubic flux is "less" nonlinear than the quadratic flux.
Notice that, with this definition,  a linear flux is not nonlinear: $d=+\infty$
with the natural convention $\min(\emptyset)=+\infty$.

 This definition is equivalent  to  Definition \ref{defLPT} for $C^\infty$ flux 
with $  \alpha = \dfrac{1}{d}$  (\cite{BJ,Ju4}). 
Therefore  the Lions-Perthame-Tadmor  parameter $\alpha$ is  for smooth flux the inverse of an integer. 
\\

 The conjectured smoothing effect \eqref{conj} is proved for the first time  in fractional $BV$ spaces 
 for the  class  of    nonlinear  (degenerate) convex   fluxes (\cite{BGJ6}). 
\begin{definition}{\bf [Nonlinear degenerate convex  flux, \cite{C,BGJ6}]} \label{defHDR}
$\mbox{}$\\
Let $f$ belong to $C^1(I,\mathbb{R})$ where $I$ is an interval of $\mathbb{R}$.  We say that the degeneracy of $f$ 
on $I$ is at least $p$ if the continuous derivative  $a(u)=f'(u)$ satisfies: 
\begin{equation}\label{indeg}
   0  < \inf_{I \times I} \frac{|a(u) - a(v)|}{|u-v|^p}
\end{equation} 
The lowest real number $p$, if there  exists, is called the degeneracy of $f$ on $I$.  If there is no $p$ such that  \eqref{indeg} is satisfied, we set $p=+ \infty$.\\
Let $f\in C^2(I)$. We say that a real number $y\in I$ is a degeneracy point of $f$ on $I$ if
$f''(y)=0$ (i.e. $y$ is a critical point of $a$).
\end{definition}
 For instance, if $f$ is the  power-law flux on $[-1,1]$: $f(u)=| u|^{1+\alpha}$ where $\alpha >0$,  then  the   degeneracy is 
$p = \max(1,\alpha)$,  (\cite{C,BGJ6}).
\begin{remark}
 Definition \ref{defHDR}    implies    the convexity (or the concavity) of the flux $f$. 
\end{remark}
Indeed, by definition there exists $C>0$ such that $|f'(u) - f'(v)| \geq C |u-v|^p$.  
Hence  the difference $ f'(u) - f'(v)$  never  vanishes for $u\neq v$. 
Since the flux is continuous, it has got a constant sign   for $u>v$, 
which implies the monotonicity of $f'$ and then the convexity  (or the concavity) of the  flux. 

\begin{remark}
 Definition \ref{defHDR} is less general than Definition \ref{defLPT}.  
 Nevertheless,  if  $f$ satisfies (\ref{indeg})  then it also satisfies   
  (\ref{nfLPT})  with  $\alpha = \dfrac{1}{p}$ ,     and also  (\ref{eqJ4}) with $d=p$  when  $f$ is smooth. 
  \end{remark}

 The paper is organized as follows. 
The sequence  given in Section \ref{sHFW}  is  exactly uniformly bounded in the Sobolev space conjectured in \cite{LPT}. 
 Furthermore, this sequence is  unbounded in all smoother Sobolev spaces. 
In  Section  \ref{sos},     we build   solutions
 with    the  suitable  regularity \eqref{conj}.



\section{Supercritical geometric optics} \label{sHFW}
\smallskip
$\mbox{}$\\

We give   a sequence of high frequency waves with small amplitude exactly uniformly bounded in the Sobolev space conjectured in \cite{LPT}. 
The construction uses a WKB  expansion (\cite{CJR,Rauch}).
\begin{theorem}
Let $f \in C^\infty(K,\mathbb{R})$ be a nonlinear flux with degeneracy $d$  defined by (\ref{eqJ4}). 
There exists a constant state $\underline{u} \in  K$  
such that
for any smooth periodic   function $U_0$ 
satisfying for all $0 < \varepsilon \leq 1$, for all $x \in \mathbb{R}$,
$u_0^\varepsilon(x)= \underline{u} + \varepsilon U_0 \left(\frac{x}{\varepsilon^d}\right) \in K$,
the following properties hold:
\begin{enumerate}
\item  there exists  a positive time  $T$ such that
 the entropy solution  $u^\varepsilon $ of equation (\ref{eqcl}) with $u_0=u_0^\varepsilon$
is  smooth on $[0,T]\times \mathbb{R}$  for all $0 < \varepsilon \leq 1$ ,
\item 
the sequence $(u^\varepsilon)$ is uniformly  bounded in $W_{loc}^{s,1}([0,T]\times \mathbb{R})$  for $s= \alpha= \dfrac{1}{d}$ 
and unbounded for $s > \alpha$  when  $U_0' \neq 0 $ a.e. 
\end{enumerate}
\end{theorem}
The  key point is to construct a  sequence of very   high frequency waves near the  state $\underline{u}$ 
where the maximum in \eqref{eqJ4} 
is reached. Next we compute  the  optimal Sobolev bounds  uniformly with respect to $\varepsilon$   on the  WKB  expansion:
$$ u^\varepsilon(t,x)=  \underline{u} + \varepsilon \;  U\left( t, \dfrac{\varphi(t,x)}{\varepsilon^d}\right) + \varepsilon \;  r_\varepsilon (t,x).$$
To estimate the remainder  in Sobolev norms, we build a smooth sequence of solutions. 
It is quite surprising to have such smooth sequence  on uniform time strip $[0,T]$. 
Indeed, it is a sequence of solutions with no entropy production, without shock.  
But for any higher  frequency,  the  life span  $T_\varepsilon $  of  $u_\varepsilon$  as  a continuous solution   goes towards  $ 0$  and oscillations are canceled (\cite{Ju4}).  
Thus the construction is optimal. 
\begin{remark} The uniform life span of the smooth sequence $(u^\varepsilon)$ is at least 
$$T \sim \dfrac{1}{\sup_\theta \left| \dfrac{ d \,  U_0} {d \, \theta} \right|},$$  as one can see in \cite{Ju4}. 
So we can build such smooth sequence for  any   large time $T$ and any  non constant initial periodic profile $U_0$ small enough in $C^1$. 
But we cannot take $T=+\infty$ since shocks always occur when $U_0$ is not constant. 
\end{remark}

\begin{remark} \label{allal}   For $C^\infty$ flux,  the  parameter $\alpha$ in Definition \ref{defLPT} is always the inverse of an integer. 
 To get supercritical geometric optics  expansions  for all $\alpha \in ]0,1]$ and not only  $\alpha \in \left\{\dfrac{1}{n}, \;  n \in \mathbb{N}^* \right\}$,  
 we shall consider  power-law flux $f(u)= |u|^{1+ p}$, where  $p = \dfrac{1}{\alpha} \in [1,+\infty[$, as in \cite{DW}.
 In this  case, $\underline{u}=0$ and  the sequence  is simply 
$u^\varepsilon(t,x)= \varepsilon  U\left (t , \dfrac{x}{\varepsilon^p}\right) $, 
the exact entropy solution of (\ref{eqcl}) and  $U(t,\theta)$ is the entropy solution of   
$ \partial_t  U + \partial_\theta |U|^{1+p} =0$, $U(0,\theta)= U_0(\theta)$. 
\end{remark}

\pr  We give a sketch of the proof  (see \cite{Ju4} for more details).
\begin{itemize}
\item   Existence of $\underline{u}$: 
   the map $u \rightarrow \min\{ k \geq 1, \; f^{(1+k)} (u) \neq 0\}$ is upper semi-continuous, 
 so  it  achieves its maximum  on the compact $K$.
\item  WKB expansion (\cite{DM,Ju2,CJR,Ju4}):
 we plug the ansatz   $$u^\varepsilon(t,x)=  \underline{u} + \varepsilon U_\varepsilon \left( t, \frac{\varphi(t,x)}{\varepsilon^d}\right) $$
into (\ref{eqcl}).  
 Notice that the exact profile $U_\varepsilon$ depends on $\varepsilon$.
 \\
 Set $\lambda= f'(\underline{u})$ and  $b= \dfrac{f^{(1+d)} (\underline{u})}{(1+d)!} \neq 0$. 
After simplification, the Taylor expansion of the flux   
 $f(\underline{u} + \varepsilon U_\varepsilon)
=  f(\underline{u}) + \varepsilon \,  \lambda \,  U_\varepsilon  + \varepsilon^{1+d} \,  b \,  U_\varepsilon^{1+d}
-  \varepsilon^{2+d}\, R_\varepsilon(U_\varepsilon) $ gives an equation for the exact profile $U_\varepsilon$ and the phase $\varphi$:
 \begin{equation} \label{eqexactprofile}
\frac{\partial   U_\varepsilon}{\partial t } +  b \frac{\partial   U_\varepsilon^{1+d}}{\partial \theta }
 =  \varepsilon \frac{\partial   R_\varepsilon( U_\varepsilon)}{\partial \theta },
 \quad U_\varepsilon(0,\theta)=U_0(\theta), 
\quad \varphi(t,x)= x - \lambda t .
\end{equation}
The profile,  which does not depend on $\varepsilon$, is
 \begin{equation} \label{eqprofile}
\frac{\partial   U}{\partial t } +  b \frac{\partial   U^{1+d}}{\partial \theta } = 0, \qquad U(0,\theta)=U_0(\theta).
\end{equation}

\item    Existence of smooth solutions  for a time $T>0$ independent of $\varepsilon$: 
       it is a consequence of the method of characteristics. 
       Indeed, the characteristics of equation  (\ref{eqexactprofile}) are a small perturbation 
  of characteristics of  equation (\ref{eqprofile}).  
\item    Approximation  in $ C^1 ([0,T]\times \mathbb{R})$: it  comes again
  from  the method of characteristics since $\varepsilon R_\varepsilon \rightarrow 0$.
\\
  Notice that the expansion is valid   in $L^1_{loc}$ after shock waves  (\cite{CJR}). But it is not enough to estimate the Sobolev norms. 
\item   Sobolev estimates:    
roughly speaking,  the order of growth of the $s$ fractional derivative  $   \dfrac{d^s}{d x^s}  U_0 \left(\frac{x}{\varepsilon^d}\right) $
 is  $ \varepsilon^{ -sd} .$ 
For the profile $U$, this estimate is propagated along the characteristics on $[0,T]$. 
 We have the same estimate for $U_\varepsilon$ since $U_\varepsilon$ is near $U$ in $C^1$. 
Then we get  the Sobolev bounds for $u_\varepsilon$.
\end{itemize}
\cqfd

\section{Oscillating solutions} \label{sos} 
\smallskip
$\mbox{}$\\

In this section we give exact continuous solutions with   the  Sobolev regularity conjectured in \cite{LPT}. 
Indeed, we choose a  suitable initial data  such that the regularity  is not spoiled by the nonlinearity of the flux for a positive time $T$.
Furthermore,   
the conjectured smoothing effect is proved for the first time  in fractional $BV$ spaces (\cite{BGJ6})
 for the   degenerate convex  class of nonlinear flux given by Definition \ref{defHDR}. 
The next theorem shows  the optimality of this smoothing effect. 
  The optimality was also given in \cite{DW} in  Besov spaces framework.
Let us introduce the  $BV^s$ spaces. 
\begin{definition}[Fractional $BV$ spaces]
$\mbox{}$\\
Let I be a non empty interval of $\mathbb{R}$. A partition  $\sigma$ of the interval $I$  is a finite ordered subset: 
 $\sigma=\{x_0, x_1, \cdots,x_n\} \subset I$, $x_0 < x_1 < \cdots < x_n$.
 We denote by $S(I)$ the set of all  partitions of $I$. 
 Let $s$ belong to $]0,1]$ and $ p =\dfrac{1}{s} \geq 1$.  The s-total variation  of a real function  $u$ on $I$ is 
\begin{eqnarray*}
   TV^s u\{I\} & = & \displaystyle{  \sup_{\sigma \in S(I)} \sum_{k=1}^n \left|u(x_k) - u(x_{k-1})\right|^p. }
\end{eqnarray*}
$BV^s(I)$ is the space of real functions $u$  such that $TV^s u \{I\} < + \infty$.
\end{definition}
 $BV^s$ spaces are introduced in   \cite{BGJ6} for applications to  conservation laws.  
 These spaces  measure the regularity of regulated functions: $BV=BV^1 \subset BV^s \subset L^\infty$. 
 Indeed, $BV^s(K)$ is very close to the Sobolev space   $W^{s,1/s}(K)$ (\cite{BGJ6}): 
\begin{itemize}   
\item   $BV^s(K) \subset W^{s - \eta, 1/s}(K)$  for all $0 < \eta < s $ . 
\item $BV^s(K) \neq W^{s,1/s}(K) $
\end{itemize}

We now give continuous functions which have  the $BV^s$ regularity.
\begin{proposition}[A continuous $BV^s$ function \cite{C}]\label{pBVsex} 
$\mbox{}$\\
Let $ 0 < s < 1$, $0 < \eta < 1 -s $ and let  $g=g_{s,\eta}$ be the real function defined on $[0,1]$ by 
 $g(0)=0$ and
for all $x \in ]0,1]$ :
\[
{\displaystyle g(x)=x^{b}\cos\left(\dfrac{\pi}{x^{c}}\right)}, \quad \textrm{where}\;\; b= s + \dfrac{s^2}{\eta}\;\;\textrm{and}\;\ c=\dfrac{s}{\eta}.
\]
The function $g$ belongs to $ BV^{s}([0,1]) \cap C^0([0,1])$
but not to  $ BV^{s+\eta}([0,1])$.
\end{proposition}
 \noindent Notice that such example do not provide a function which belongs to $BV^s$ but not to  $ \bigcup_{\eta >0} BV^{s+\eta}$.
\pr
 The extrema of $g$ are achieved on $x_k= k^{- 1/c}$.    
 Let $p =\dfrac{1}{s} > 1$ ,  $q \leq p $ and 
$$  V_q = \sum_{k=1}^{+\infty} |g(x_{k+1}) - g(x_{k})|^q .$$ 
 Since  
$ q b / c = q  (s+\eta)$, the asymptotic behavior $ |g(x_{k+1}) - g(x_{k})|^q \sim  2^q  k^{- q b/c}$ 
 when $k \rightarrow + \infty$ 
 yields   $V_q=+ \infty $ when $q=1/(s+\eta) $ and $ V_p < + \infty $. 
First  this implies  $g \notin  BV^{s+\eta}$. 
Second,  for  such  oscillating function with diminishing amplitudes, 
we choose the optimal infinite partition  to compute the 
 s-total variation (see Proposition 2.3. p. 6 in \cite{BGJ6}). Then $g$ belongs to $BV^s$.
\cqfd


We are now able to find oscillating initial data with the critical Sobolev exponent propagated by the nonlinear conservation law (\ref{eqcl}).

\begin{theorem} \label{thCheng}
  Assume  $f\in C^{\infty}(K,\mathbb{R})$  be nonlinear in the sense of
  Definition  \ref{defJ4}. 
We denote by $d$ its degeneracy and $s =\dfrac{1}{d}$. 
 For any    $\eta >0 $ and any  time $T>0$ there exists  a  solution $u \in C^0([0,T]\times \mathbb{R},\mathbb{R})$  such that
for all $t\in[0\,,\, T]$ \\
$$ u(t,\cdot)\in BV^{s}(\mathbb{R}, \mathbb{R})  \mbox{  and } u(t,\cdot)\notin BV^{s+\eta}(\mathbb{R}, \mathbb{R}).$$
\end{theorem}

The idea follows the K-S Cheng construction (\cite{Cheng83}) with the function $g$ given in Proposition \ref{pBVsex}.

\pr
Let  $\underline{u}\in K$ a point where the maximum of degeneracy of $f$ is achieved. We
also suppose that $\underline{u}\in\overset{\circ}{K}$ 
(the proof of
Theorem \ref{thCheng}  is quite similar
if $\underline{u}\in\partial K$).
\\
We define the initial condition
$u_{0}$ by :
\[
\left\{ \begin{array}{rl}
u_{0}(x)=\underline{u} & \textrm{ if }x<0\\
\\
u_{0}(x)=\underline{u}+\delta g(x) & \textrm{ if }0\leq x\leq1\\
\\
u_{0}(x)=\underline{u}-\delta & \textrm{ if }1<x
\end{array}\right.,
\]
where $\delta>0$ is chosen such that for all $x\in[0,1]$, $\underline{u}+\delta g(x)\in K$.
Notice that for all $x\in[0,1]$, $-1\leq g(x)\leq1$ and $g(1)=-1$.

Then, following the method of characteristics, we define the function
$u(t,x)$ by :
\[
\left\{ \begin{array}{rl}
u(t,x)=0 & \textrm{ if }x<0\\
\\
u(t,x)=\underline{u}+\delta g(y) & \textrm{ if }x=y+ta(\underline{u}+\delta g(y)),\quad0\leq y\leq1\\
\\
u(t,x)=\underline{u}-\delta & \textrm{ if }1+ta(\underline{u}-\delta)<x
\end{array}\right..
\]

 $u_{0}\in BV^{s}([0\,,\,1])$ and $u_{0}\notin BV^{s+\eta}([0\,,\,1])$.
Let be $t>0$ and for all $y$, $$ \theta_{t}(y)=y+ta(\underline{u}+\delta g(y)) .$$

Considering the change of variable $y= x - a(\underline{u})t$,
we can assume  without loss  of generality  that $f'(\underline{u})=a(\underline{u})=0$.
Since $f\in C^{\infty}(K,\mathbb{R})$, we derive from a Taylor expansion that
\[
a(u)=\frac{1}{d!}\left(a^{(d)}(\underline{u})(u-\underline{u})^{d}+\int_{\underline{u}}^{u}(u-s)^{d}a^{(1+d)}(s)ds\right).
\]
Defining
\[
I_{n}(y)=\frac{1}{d!}\int_{0}^{1}(1-r)^{d}a^{(1+d)}(\underline{u}+r\, \delta \,  g(y))dr,\] 
\[ J_{n}(y)=\frac{1}{d!}\int_{0}^{1}r(1-r)^{d}a^{(2+d)}(\underline{u}+r\, \delta \,  g(y))dr,
\]
we get then :
\[
\theta_{t}(y)=y+t\delta^{d}g(y)^{d}\left(\frac{1}{d!}a^{(d)}(\underline{u})+\delta g(y)I_{n}(y)\right).
\]
Note that $g$, $I_{n}$, $J_{n}$ are bounded on $[0,1]$.
\\  For $y\neq0$, since  $b\, d=1+c$, we have
${\displaystyle \frac{\left|g(y)\right|^{d}}{y}=O\left(y^c \right)}$
at $0$.
Thus $\theta_{t}$ is differentiable at $0$ and $\dfrac{d\theta_{t}}{dy}(0)=1$.
For $y\neq0$, we have 
 $$\dfrac{d\theta_{t}}{dy}(y)=1+t\delta^{d}h_{n}(y),$$
where
\[
h_{n}(y)=g(y)^{d-1}g'(y)\left(\frac{1}{(d-1)!}a^{(d)}(\underline{u})+(d+1)\delta g(y)I_{n}(y)+\delta^{2}g(y)^{2}J_{n}(y)\right).
\]
For $y\neq0$, since  $b\, d=1+c$,  we have 
\begin{eqnarray*}
{\displaystyle g(y)^{d-1}g'(y)} &= & { \left( y^{b_{}}\cos\left(\dfrac{\pi}{y^{c_{}}}\right) \right)^{d-1} \left(b_{} \, y^{b_{}-1}\cos\left(\dfrac{\pi}{y^{c_{}}}\right)+\pi c_{}y^{b-c_{}-1} \sin\left(\dfrac{\pi}{y^{c_{}}}\right)\right)},
\\
{\displaystyle \left|g(y)^{d-1}g'(y)\right| } & \leq  & {  \left|\cos\left(\dfrac{\pi}{y^{c_{}}}\right)\right|^{d-1}\left(b_{}\left|y\right|^{c}\left|\cos\left(\dfrac{\pi}{y^{c_{}}}\right)\right|+\pi c_{}\left|\sin\left(\dfrac{\pi}{y^{c_{}}}\right)\right|\right)}.
\end{eqnarray*}
Thus $g(y)^{d-1}g'(y)$ is bounded on $[0,1]$.

As $h_{n}$ is bounded on $[0,1]$, there exists $T_\delta>0$
such that for all $y\in[0,1]$ and for all $t\in]0,T]$, $\dfrac{d\theta_{t}}{dy}(y)>0$. 
Notice that $ \displaystyle{\lim_{\delta \rightarrow 0} T_\delta = + \infty}$.
We can take $\delta >0$ small enough such that $T_\delta > T$.

Thus for all $t\in]0,T]$, $\theta_{t}$ is an homeomorphism between
$[0,1]$ and $[0,1+ta(\underline{u}-\delta)]$. Then $u(t,x)$ is a continuous  solution
of equation (\ref{eqcl}) on $[0,T] \times \mathbb{R} $.
Furthermore, since $u_{0}\in BV^{s}(I)$ and $u_{0}\notin BV^{s+\eta}(I)$,
where $I=[0,1]$, we deduce that for all $t\in]0,T]$, $u(t,\cdot)\in BV^{s}(J)$
and $u(t,\cdot)\notin BV^{s+\eta}(J)$, where $J=\theta_{t}(I)=[0,1+ta(\underline{u}-\delta)]$.
Finally, as $u(t,\cdot)$ is constant outside $J$, we have proved
that $u(t,\cdot)\in BV^{s}(\mathbb{R})$ and $u(t,\cdot)\notin BV^{s+\eta}(\mathbb{R})$.
\cqfd

\begin{remark}      
As in Remark \ref{allal}, Theorem \ref{thCheng} is restricted  for critical exponent $s$ such that $\dfrac{1}{s} \in \mathbb{N}$. 
To obtain all exponent $s \in ]0,1]$, following \cite{DW},
we can  consider  a power-law  flux with $p=\dfrac{1}{s}$: 
$
f(u)=|u|^{1+p}.
$
Our construction is quite similar as in the proof of Theorem  \ref{thCheng} with $\underline{u}=0$ and $\delta > 0 $ small enough.
\end{remark}



     
\end{document}